\documentclass[titlepage,twoside,12pt]{article}
\usepackage{amssymb}
\usepackage{amsfonts}
\begin{document}

{\bf \Large An extension of the Eilenberg - Mac Lane \\ \\ concept of category in the case
of nonrigid \\ \\ structures} \\ \\

{\bf Elem\'{e}r E Rosinger} \\
Department of Mathematics \\
and Applied Mathematics \\
University of Pretoria \\
Pretoria \\
0002 South Africa \\
eerosinger@hotmail.com \\ \\

{\bf Abstract} \\

Nonrigid mathematical structures may no longer form usual Eilenberg - Mac Lane categories,
but more general ones, as illustrated by pseudo-topologies. \\
A rather general concept of pseudo-topology was used in constructing differential algebras of
generalized functions containing the Schwartz distributions, [4-6,8-11]. These algebras proved
to be convenient in solving large classes of nonlinear partial differential equations, see
[12,13] and the literature cited there, as well as section 46F30 in the Subject Classification
2000 of the American Mathematical Society, at www.ams.org/index/msc/46Fxx.html \\
The totality of such pseudo-topologies no longer constitutes a usual Eilenberg - Mac Lane
category, but an extended one which is presented here. \\
Other nonrigid mathematical structures are mentioned and treated briefly. \\
This is a revised and augmented version of the earlier published paper [7]. \\ \\

\newpage

{\large \bf 1. Nonrigid structures} \\

Instead of a general definition, we shall give an example of what we mean by a {\it nonrigid}
mathematical structure. \\
Further examples of nonrigid structures can be seen in section 5. \\

With the usual Hausdorff-Kuratowski- Bourbaki, or in short HKB, concept of topology, [1,14],
all the topological entities, such as open sets, closed sets, compact sets, convergence,
continuous functions, etc., are {\it uniquely} determined by {\it one single} such concept,
for instance, that of open sets. In other words, if one for instance chooses to start with the
open sets, then all the other topological entities can be defined in a unique manner based on
the given open sets, [5, p. 225], [14]. \\

In contradistinction to that usual situation, in pseudo-topological structures there is a {\it
relative independence} between various topological entities. Namely

\begin{itemize}

\item the connection between topological entities is no longer rigid to the extreme as in the
case of the usual HKB concept of topology, where one of the topological entities determines in
a unique manner all the other ones; instead, the various topological entities are only
required to satisfy certain compatibility conditions, [5, p. 225], [14],

\item there exist, however, several topological entities of prime importance, among them, that
of {\it Cauchy-ness} which is a binary relation between two arbitrary sequences, nets, filters,
etc., and it was first introduced in [5], see also [6,8-11,14].

\end{itemize}

And it is precisely because of that nonrigid character that the totality of pseudo-topologies
does {\it no longer} constitute a usual Eilenberg - Mac Lane category, but one which is more
general, namely, as given in Definition 1, next.  \\ \\

\newpage

{\large \bf 2. An extended concept of category} \\

As often in Category Theory, [3], we shall use the concept of {\it class} which, intuitively,
is supposed to extend that of set. \\

{\bf Definition 1} \\

Let $\Phi$ be a class of objects {\bf F}. Each such object {\bf F} is a set of elements
${\cal F}$. Further, each ${\cal F}$ is a set of elements $f$ called mappings. \\
Let now $\prod$ be a class of pairs ({\bf F}, {\bf G}), where {\bf F}, {\bf G} are objects in
$\Phi$. \\
Finally, let $\circ$ be an operation of composition of suitable pairs of mappings $f, g$. \\

We say that $( \Phi, \prod, \circ )$ is an {\it extended category}, if and only if the
following conditions are satisfied : \\

(C1) $~~~ \begin{array}{l}
             \forall~~ ( {\bf F}_1, {\bf F}_2 ) \in \prod ~: \\ \\
             \exists~~ {\bf F}_3 \in \Phi ~: \\ \\
             \forall~~ {\cal F}_1 \in {\bf F}_1,~ {\cal F}_2 \in {\bf F}_2 ~: \\ \\
             \exists~~ {\cal F}_3 \in {\bf F}_3 ~: \\ \\
             \forall~~ f_1 \in {\cal F}_1,~ f_2 \in {\cal F}_2 ~: \\ \\
             ~~~~ f_1 \circ f_2 \in {\cal F}_3
          \end{array} $ \\ \\

also \\

\newpage

(C2) $~~~ \begin{array}{l}
             \forall~~ ( {\bf F}_1, {\bf F}_2 ), ( {\bf F}_2, {\bf F}_3 )\in \prod ~: \\ \\
             \exists~~ {\bf F}_4 \in \Psi ~: \\ \\
             \forall~~ {\cal F}_1 \in {\bf F}_1,~ {\cal F}_2 \in {\bf F}_2,~
                                                      {\cal F}_3 \in {\bf F}_3 ~: \\ \\
             \exists~~ {\cal F}_4 \in {\bf F}_4 ~: \\ \\
             \forall~~ f_1 \in {\cal F}_1,~ f_2 \in {\cal F}_2,~ f_3 \in {\cal F}_3 ~: \\ \\
              ~~~~ f_1 \circ ( f_2 \circ f_3 ) ~=~ ( f_1 \circ f_2 ) \circ f_3
          \end{array} $ \\ \\

as well as \\

(C3) $~~~ \begin{array}{l}
              \forall~~ {\bf F} \in \Phi ~: \\ \\
              \exists~~ {\bf F}_0,~ {\bf F}^0 \in \Phi ~: \\ \\
              ~~~~ (C3.1)~~~ ( {\bf F}_0, {\bf F} ),~ ( {\bf F}, {\bf F}^0 ) \in \prod ~: \\ \\
              ~~~~ (C3.2)~~~ \exists~~ f_0 \in {\cal F}_0 \in {\bf F}_0,~
                                                 f^0 \in {\cal F}^0 \in {\bf F}^0 ~: \\ \\
              ~~~~~~~~~~~~~~~\forall~~ f \in {\cal F} \in {\bf F} ~: \\ \\
              ~~~~~~~~~~~~~~~~~~~ f_0 \circ f ~=~ f \circ f^0 ~=~ f
           \end{array} $

\hfill $\Box$ \\

The meaning of this definition is clarified in the next section, when it is shown how it
contains as a particular case the usual notion of Eilenberg - Mac Lane category. \\ \\

{\large \bf 3. The particular case of the Eilenberg - Mac Lane category} \\

It is known, [2,3], that a usual Eilenberg - Mac Lane category, [2], can be given as a pair
$( {\cal F}, \circ )$, where ${\cal F}$ is a class of mappings, while $\circ$ is the
composition of mappings. In such a case one considers the class ${\cal O}$ of objects given by
the identity mappings in ${\cal F}$. Further, for any two objects $A,~ B \in {\cal O}$, one
denotes by $Hom ( A, B )$ the set of mappings from $A$ to $B$. Then the composition $\circ$ of
mappings operates according to $Hom ( A, B ) \times Hom ( B, C ) \ni ( f, g ) \longmapsto g
\circ f \in Hom ( A, C )$. \\

We show now how to the usual Eilenberg - Mac Lane category $( {\cal F}, \circ )$ one can
associate an {\it extended category} in the sense of Definition 1. \\
As mentioned, this result also helps in clarifying the meaning of Definition 1, in view of the
simple manner in which a usual Eilenberg - Mac Lane category proves to be an extended category
in the sense of Definition 1 above. \\

Indeed, let us consider

\begin{itemize}

\item $\Phi$ being the class of all sets with one single element ${\bf F} = \{~ Hom ( A, B )
~\}$, where $A,~ B \in {\cal O}$,

\item $\prod$ being the class of all pairs $(~ {\bf F}_1, {\bf F}_2 ) = ( \{~ Hom ( A, B )
~\}, \{~ Hom ( B, C ) ~\} ~)$, where $A, B, C \in {\cal O}$.

\end{itemize}

We note the particular nature of the above objects ${\bf F}$, each of which is a set with one
single element given by the set $Hom ( A, B )$. This is a consequence of the fact that the
usual Eilenberg - Mac Lane categories do not model nonrigid mathematical structures. Indeed,
as seen in section 4, when we deal with pseudo-topologies, which are nonrigid mathematical
structures, the corresponding objects ${\bf F} ( \Sigma, \Sigma^\prime )$ will typically be
rather large sets, and thus no longer having only one single element. \\

Now we obtain by a step by step direct verification, see [7, pp.882] \\

{\bf Theorem 1} \\

$( \Phi, \prod, \circ )$ is an {\it extended category} in the sense of Definition 1. \\

{\bf Proof} \\

We shall show that $( \Phi, \prod, \circ )$ satisfies the above conditions (C1) - (C3) in
Definition 1. \\

For condition (C1), let us take any $( {\bf F}_1, {\bf F}_2 ) \in \prod$. Then clearly \\

$\exists~~ A, B, C \in {\cal O} ~:~ {\bf F}_1 = \{~ Hom ( A, B ) ~\},~ {\bf F}_2 =
\{~ Hom ( B, C ) ~\}$ \\

Let then ${\bf F}_3 = \{~ Hom ( A, C ) ~\}$, thus ${\bf F}_3 \in \Phi$, and it satisfies (C1).
Indeed, let ${\cal F}_1 \in {\bf F}_1,~ {\cal F}_2 \in {\bf F}_2$. Then we must have
${\cal F}_1 = Hom ( A, B ),~ {\cal F}_2 = Hom ( B, C )$. If now we take ${\cal F}_3 = Hom ( A,
C )$, then obviously ${\cal F}_3 \in {\bf F}_3$, and \\

$\forall~~ f_1 \in {\cal F}_1,~ f_2 \in {\cal F}_2 ~:~ f_1 \circ f_2 \in {\cal F}_3$ \\

To show (C2), we take $( {\bf F}_1, {\bf F}_2 ),~( {\bf F}_2, {\bf F}_3 ) \in \prod$. Then
obviously \\

$~~~ \begin{array}{l}
        \exists~~ A, B, C, D \in {\cal O} ~: \\ \\
        ~~~~ {\bf F}_1 = \{~ Hom ( A, B ) ~\},~ {\bf F}_2 = \{~ Hom ( B, C ) ~\},~
        {\bf F}_3 = \{~ Hom ( C, D ) ~\}
     \end{array} $ \\

Let ${\bf F}_4 = \{~ Hom ( A, D ) ~\}$, then clearly ${\bf F}_4 \in \Phi$, and (C2) is
satisfied. Indeed, let ${\cal F}_1 \in {\bf F}_1,~ {\cal F}_2 \in {\bf F}_2,~ {\cal F}_3 \in
{\bf F}_3$, then ${\cal F}_1 = Hom ( A, B ),~ {\cal F}_2 = Hom ( B, C ),~ {\cal F}_3 = Hom
( C, D )$. Let now ${\cal F}_4 = Hom ( A, D )$, then ${\cal F}_4 \in {\bf F}_4$, and clearly \\

$~~~ \begin{array}{l}
         \forall~~ f_1 \in {\cal F}_1,~ f_2 \in {\cal F}_2,~ f_3 \in {\cal F}_3 ~: \\ \\
         ~~~~ f_1 \circ ( f_2 \circ f_3 ) ~=~ ( f_1 \circ f_2 ) \circ f_3
     \end{array} $ \\

At last, we show that (C3) is also satisfied. Let ${\bf F} \in \Phi$. Then \\

$\exists~~ A, B \in {\cal O} ~:~ {\bf F} ~=~ \{~ Hom ( A, B ) ~\} $ \\

Now we can take ${\bf F}_0 = \{~ Hom ( A, A ) ~\},~ {\bf F}^0 = \{~ Hom ( B, B ) ~\} \in \Phi$.
Then clearly $( {\bf F}_0, {\bf F} ),~ ( {\bf F}, {\bf F}^0 ) \in \prod$. Let ${\cal F}_0 =
Hom ( A, A ),~ {\cal F}^0 = Hom ( B, B )$, then ${\cal F}_0 \in {\bf F}_0,~ {\cal F}^0 \in
{\bf F}^0$. Finally, let $id_A \in Hom ( A, A ),~ id_B \in Hom ( B, B )$ the respective
identity mappings. Denoting now $f_0 = id_A,~ f^0 = id_B$, it follows that $f_0 \in {\cal
F}_0,~ f^0 \in {\cal F}^0$, and thus (C3) is satisfied. Indeed, if ${\cal F} \in {\bf F}$,
then ${\cal F} = Hom ( A, B )$, hence \\

$\forall~~ f \in {\cal F} ~:~ f_0 \circ f ~=~ f \circ f^0 ~=~ f$ \\ \\

{\large \bf 4. The extended category of pseudo-topologies} \\

The concept of extended category in the sense of Definition 1 was motivated by the
pseudo-topologies introduced and used in [4-6], see also [8-11,14]. Here we shortly recall the
needed details. \\

Let $\Sigma, \Sigma^\prime$ be pseudo-topologies on the nonvoid sets $E$ and $E^\prime$,
respectively. Let $\Psi$ and $\Psi_u$ be the sets of mappings from $E$ to $E^\prime$ which are
continuous, respectively uniformly continuous, from the pseudo-topology $\Sigma$ on $E$ to the
pseudo-topology $\Sigma^\prime$ on $E^\prime$. Then \\

$~~~ \Theta ~=~ ( E, E^\prime, \Sigma, \Sigma^\prime, \Psi, \Psi_u ) $ \\

is called a pseudo-topology on $( E, E^\prime )$, [5, p. 227]. \\

The essential point, mentioned in section 1, is that - owing to the lesser conceptual rigidity
involved - the totality of pseudo-topologies does {\it no longer} constitute a usual Eilenberg
- Mac Lane category, but one which is more general, namely as given in Definition 1 above. \\

Among others, this lesser rigidity is manifested as follows. In a pseudo-topology $\Theta =
( E, E^\prime, \Sigma, \Sigma^\prime, \Psi, \Psi_u )$ on a given pair of spaces $( E,
E^\prime )$ as above, the sets of continuous, respectively uniform continuous functions $\Psi$
and $\Psi_u$, need {\it not} be uniquely determined by $\Sigma, \Sigma^\prime$. \\
And as mentioned, it is precisely because of that fact that the totality of pseudo-topologies
does {\it no longer} constitute a usual Eilenberg - Mac Lane category, but one which is more
general, namely as given in Definition 1 above. \\

Let us now show how the class of such pseudo-topologies $\Theta = ( E, E^\prime, \Sigma,
\Sigma^\prime, \Psi, \Psi_u )$ constitutes an extended category in the sense of Definition 1
above. For that, let us consider

\begin{itemize}

\item ${\cal F} ( \Theta ) ~=~ \Phi \times \Psi_u$

\item ${\bf F} ( \Sigma, \Sigma^\prime ) ~=~ \{~ {\cal F} ( \Theta ) ~~|~~ \Theta ~=~
( E, E^\prime, \Sigma, \Sigma^\prime, \Psi, \Psi_u ) ~ \}$, where $E, E^\prime, \Sigma$ and
$\Sigma^\prime$ are fixed,

\item $\Phi^{PT}$ being the class of all ${\bf F} ( \Sigma, \Sigma^\prime )$, where $E,
E^\prime, \Sigma$ and $\Sigma^\prime$ are arbitrary,

\item $\prod^{PT}$ being the class of all pairs $( {\bf F} ( \Sigma, \Sigma^\prime_1 ),
{\bf F} ( \Sigma^\prime_2, \Sigma^{\prime \prime} ) )$, where the pseudo-topologies
$\Sigma^\prime_1, \Sigma^\prime_2$ on $E^\prime$ are in the relation $\Sigma^\prime_2 \leq
\Sigma^\prime_1$, that is, $\Sigma^\prime_2$ is more fine than $\Sigma^\prime_1$, [4],

\item the composition $\circ$ is defined by $( f, f_u ) \circ ( f^\prime, f^\prime_u ) =
( f^\prime \circ f, f^\prime_u \circ f_u )$, where the last two $\circ$ are the usual
composition of mappings.

\end{itemize}

Clearly, the above objects ${\bf F} ( \Sigma, \Sigma^\prime )$ are typically rather large sets,
and thus contain not only one single element. This, as noted, is due to the fact that pseudo-
topologies are nonrigid mathematical structures. \\

Now we obtain, [7, p. 884] \\

{\bf Theorem 2} \\

$( \Phi^{PT}, \prod^{PT}, \circ )$ is an {\it extended category} in the sense of
Definition 1. \\

{\bf Proof} ( see [7, pp. 884, 885] ). \\ \\

{\large \bf 5. A few comments} \\

One can note that even in Algebra there are nonrigid mathematical structures. For instance,
let $( R, +, . )$ be a ring. Then in principle, neither the addition "+" determines the
multiplication ".", nor multiplication determines the addition. Instead, they are relatively
independent, and only satisfy the usual compatibility conditions. \\

In the case of the usual ring $\mathbb{Z}$ of integers, multiplication is in fact determined
by addition as a repeated addition, since for every $m,\, n \in \mathbb{Z},~ m \geq 1$, we
have $m \,.\, n = n \,+\, .~.~. \,+ n$, where the respective sum has $m$ terms. \\
Thus as a ring, $\mathbb{Z}$ has a rigid mathematical structure. \\

As for arbitrary rings, at least in principle, they need not be rigid, although there is not
enough clarity about the existence of nonrigid rings in general. \\

Spaces $( \Omega, {\cal M}, \mu )$ with measure, where $\Omega$ is the underlying set, ${\cal
M}$ is a $\sigma$-algebra on it, and $\mu : {\cal M} \longrightarrow \mathbb{R}$ is a $\sigma$-
additive measure, are further examples of nonrigid structures, since for a given $( \Omega,
{\cal M} )$, there can in general be infinitely many associated $\mu$. \\

Topological groups, or even topological vector spaces, are typically nonrigid structures.
Indeed, on an arbitrary group, or even vector space, there may in general be many compatible
topologies, and even Hausdorff topologies. \\

Let us end by sketching briefly how topological groups can be organized so as to form an
extended category in the sense of Definition 1 above. \\

We denote a topological group by $( G, *, \tau )$, where $G$ is the underlying set, $*$ is the
group operation and $\tau$ is the compatible topology. \\
Given two topological groups $( G, *, \tau ),~ ( G^\prime, *^\prime, \tau^\prime )$, we denote
by ~${\cal C} ( ( G, *, \tau ), ( G^\prime, *^\prime, \tau^\prime ) )$ the set of continuous
group homomorphisms $f : G \longrightarrow G^\prime$. \\
Now we consider

\begin{itemize}

\item the objects \\ ${\bf F} ( ( G, * ), ( G^\prime, *^\prime ) ) ~=~
      \left \{~ {\cal C} ( ( G, *, \tau ), ( G^\prime, *^\prime, \tau^\prime ) ) ~~
         \begin{array}{|l}
            ~\tau ~\mbox{and}~ \tau^\prime ~\mbox{are compatible} \\
            ~\mbox{topologies with}~ ( G, * )~ \mbox{and} \\
            ~( G^\prime, *^\prime ), ~\mbox{respectively}
         \end{array} ~\right \}$

\item $\Phi^G$ the class of all objects ${\bf F} ( ( G, * ), ( G^\prime, *^\prime ) )$

\item $\prod^G$ the class of all pairs of objects $( {\bf F} ( ( G, * ), ( G^\prime,
*^\prime ) ), {\bf F} ( ( G^\prime, *^\prime ), ( G^{\prime \prime}, *^{\prime \prime} ) ) )$

\item $\circ$ the usual composition of mappings

\end{itemize}

Then it can be shown that $( \Phi^G, \prod^G, \circ )$ is an extended category in the sense of
Definition 1 above,  while it is not a usual Eilenberg - Mac Lane one, since the objects
${\bf F} ( ( G, * ), ( G^\prime, *^\prime ) )$ are typically large sets, and not sets with
only one single element. \\ \\

{\bf Note} \\

It may be amusing to recall how the original paper [7] was reviewed in the Mathematical
Reviews, MR0182588(32\#71), soon after its publication. The reviewer, J R Isbell, was in those
times one of the leading specialists in several related fields, and among others, he wrote :
"... The paper is written in poetry with a prose introduction and footnote. The footnote
should not be overlooked." \\
Needless to say, I found that comment somewhat less than comforting, since it had in no way
been in my intention to write poetry ... \\
Several of my colleagues, however, tried to reassure me that the review was in fact favourable,
and that they themselves wished to receive a similar one ... \\
So much for ... poetry and mathematics ... \\ \\

\end{document}